\documentclass[11pt]{article}

\usepackage{latexsym,color,graphicx,amsthm,amsmath,amssymb}
\usepackage{lineno}

\def\reals{{\mathbb R}}

\def\A{{\cal A}}
\def\C{{\cal C}}

\def\eps{{\varepsilon}}
\def\polylog{{\rm polylog}}
\newtheorem{theorem}{Theorem}[section]

\makeatletter
\newcommand{\ProofEndBox}{{\ifhmode\unskip\nobreak\hfil\penalty50 \else
          \leavevmode\fi\quad\vadjust{}\nobreak\hfill$\Box$
            \finalhyphendemerits=0 \par}}
\makeatother

\newcommand{\ignore}[1]{}
\begin{document}

\begin{titlepage}


\title{Distinct and repeated distances on a surface and incidences between points and spheres \thanks{%
Work on this paper by Noam Solomon and Micha Sharir was supported by
Grant 892/13 from the Israel Science Foundation. Work by Micha
Sharir was also supported by Grant 2012/229 from the U.S.--Israel
Binational Science Foundation, by the Israeli Centers of Research
Excellence (I-CORE) program (Center No.~4/11), and by the Hermann
Minkowski-MINERVA Center for Geometry at Tel Aviv University. }}

\author{
Micha Sharir\thanks{%
School of Computer Science, Tel Aviv University,
Tel Aviv 69978, Israel.
{\sl michas@tau.ac.il} }
\and
Noam Solomon\thanks{%
School of Computer Science, Tel Aviv University,
Tel Aviv 69978, Israel.
{\sl noam.solom@gmail.com} }
}

\maketitle

\begin{abstract}
In their seminal paper from 2010, Guth and Katz~\cite{GK2} proved
that the number of distinct distances determined by a set of $n$
points in $\reals^2$ is $\Omega(n/\log n)$, thus (almost) settling
Erd{\H o}s's distinct distances problem, open for nearly 65 years. In
$\reals^3$, it is conjectured that a set of $n$ points determines at
least $\Omega(n^{2/3})$ distinct distances. This bound is best possible
as it is attained by the vertices of the $n^{1/3}\times n^{1/3} \times n^{1/3}$
integer grid. This problem however is still wide open, for many years. The
best known lower bound is due to Solymosi and Vu~\cite{SoVu}.

In this paper we show that the number of distinct distances determined
by a set of $n$ points on a constant-degree two-dimensional algebraic variety $V$
(i.e., a surface) in $\reals^3$ is at least $\Omega\left(n^{7/9}/\polylog \,n\right)$.
This bound is significantly larger than the conjectured bound
$\Omega(n^{2/3})$ for general point sets in $\reals^3$.

We also show that the number of unit distances determined by $n$
points on a surface $V$, as above, is $O(n^{4/3})$, a bound that
matches the best known planar bound, and is worst-case tight in 3-space.
This is in sharp contrast with the best known general bound $O(n^{3/2})$
for points in three dimensions.

We also obtain sharp bounds for bipartite versions of the distinct
distances and the repeated distances problems.

To prove these results, we establish an improved upper bound for the
number of incidences between a set $P$ of $m$ points and a set $S$
of $n$ spheres, of arbitrary radii, in $\reals^3$, provided that the
points lie on an algebraic surface $V$ of constant degree, which does not have
linear or spherical components. Specifically, the bound is
$$
O\left( m^{2/3}n^{2/3} + m^{1/2}n^{7/8}\log^\beta(m^4/n) + m + n +
\sum_{c} |P_{c}|\cdot |S_{c}| \right) ,
$$
where the constant of proportionality and the constant exponent
$\beta$ depend on the degree of $V$, and where the sum ranges over
all circles $c$ that are fully contained in $V$, so that, for each
such $c$, $P_c = P\cap c$ and $S_c$ is the set of the spheres of $S$
that contain $c$. In addition, $\sum_{c} |P_{c}| = O(m)$ and
$\sum_{c} |S_{c}| = O(n)$.

This bound too improves upon earlier known bounds. These have been obtained
for arbitrary point sets but only under severe restrictions about
the spheres, which are dropped in our result.

Another interesting application of our result is an incidence bound
for arbitrary points and spheres in 3-space, where we improve and
generalize the previous work of Apfelbaum and Sharir~\cite{ApS2}.
\end{abstract}

\noindent {\bf Keywords.} Combinatorial geometry, incidences, the
polynomial method, algebraic geometry, distinct distances.

\end{titlepage}

\section{Introduction}

\paragraph{Incidences between points and spheres.}
Let $P$ be a set of $m$ points, and $S$ a set of $n$ spheres of
arbitrary radii in $\reals^3$. Assume that $P$ is contained in some
two-dimensional algebraic variety (surface) $V$ of constant degree $D$,
which does not have any planar or spherical components. We wish to bound the size $I(P,S)$ of the
\emph{incidence graph} $G(P,S)$, whose edges connect all pairs
$(p,s)\in P\times S$ such that $p$ is incident to $s$. In general,
and in this special setup too, $I(P,S)$ can be as large as the
maximum possible value $|P|\cdot|S|$, by placing all the points of
$P$ on a circle, and make all the spheres of $S$ contain this
circle, in which case $G(P,S) = P\times S$. The bound that we are
going to obtain will of course acknowledge this possibility, and
will be of the form $I_0(P,S) + \sum_i |P_i|\cdot |S_i|$, where, for
each $i$, $P_i\subseteq P$, $S_i\subseteq S$, and $P_i\times S_i
\subseteq G(P,S)$. Moreover, each subgraph $P_i\times S_i$ is
induced by a circle contained in $V$ that contains all the points of
$P_i$ and is contained in all the spheres of $S_i$. Informally, the
residue term $I_0(P,S)$ bounds the number of ``accidental''
incidences, those that cannot be ``explained'' in terms of large
complete bipartite subgraphs of $G(P,S)$. The quality of the bound
will be measured by $I_0(P,S)$ and by $\sum_i \left( |P_i|+|S_i| \right)$.

\paragraph{Distinct and repeated distances in $\reals^3$.}
There are two main motivations for studying point-sphere incidences.
One involves repeated and distinct distances. For repeated
distances, one draws a sphere with the given distance as radius
around each input point, and the number of incidences between the
points and spheres is exactly twice the number of repetitions of
that distance. Applications of this kind
include~\cite{CEGSW,KMSS,Za}. For distinct distances, one draws
spheres centered at the given points and having as radii all the $t$
possible distances. The number of incidences of these $nt$ spheres
with the $n$ given points is exactly $n(n-1)$, so an upper bound on
point-sphere incidences will yield a lower bound on $t$; see,
e.g.,~\cite{ANP+} and~\cite{Sze} for applications of this approach.

The point-sphere incidence approach is an effective tool for
providing lower bounds on the number of distinct distances on a
surface in three dimensions, as demonstrated below in
Theorem~\ref{dd3}.

A second, closely related class of applications involves the number
of repetitions of more involved patterns, typically congruent and
similar simplices in a given point set; see~\cite{AAPS,AgS,Bra} for
examples of such applications. These applications are not discussed
in the present paper.

\paragraph{Background.}
Earlier works on point-sphere incidences have considered the general
setup, where the points of $P$ are arbitrarily placed in $\reals^3$.
Initial partial results go back to Chung~\cite{Chung} and to
Clarkson et al.~\cite{CEGSW}, and continue with the work of Aronov
et al.~\cite{APST}.  Later, Agarwal et al.~\cite{AAPS} have bounded
the number of \emph{non-degenerate} spheres with respect to a given
point set, which was then improved by Apfelbaum and
Sharir~\cite{ApS2}.\footnote{%
  Given a finite point set $P\subset\reals^3$ and a constant $0<\eta<1$,
  a sphere $\sigma \subset \reals^3$ is called $\eta$-degenerate (with
  respect to $P$), if there exists a circle $c \subset \sigma$ such that
  $|c\cap P| \ge \eta|\sigma \cap P|$.}
Recently, Zahl~\cite{Za} gave a bound for the number of incidences
between $m$ points and $n$ spheres in $\reals^3$, when every triple
of spheres intersect at a finite set of points (which is the general
case), as part of a more general bound on the number of incidences
between points and bounded-degree surfaces in $\reals^3$ satisfying
certain favorable conditions. Zahl's bound for spheres is
$O(m^{3/4}n^{3/4}+m+n)$. (This bound was later generalized by Basu
and Sombra~\cite{BS14} to incidences between points and bounded
degree hypersurfaces in $\reals^4$ satisfying certain analogous
conditions.) The case of incidences with unit spheres have also been
studied in Kaplan et al.~\cite{KMSS}, with the same upper bound.
Other mildly related recent works include~\cite{BS,CILRD,PTV}.

The study in this paper continues a similar recent study by the authors~\cite{SS16},
involving incidences between points on a variety and planes in three dimensions.

\paragraph{Main Theorem.}
As we show in this paper, the bound can be
substantially improved when all the points of $P$ lie on a
constant-degree surface $V$. Our main result is the following
theorem.
\begin{theorem} \label{main}
Let $P$ be a set of $m$ points on some algebraic surface $V$ of
constant degree $D$ in $\reals^3$, which has no linear or spherical
components, and let $S$ be a set of $n$ spheres, of arbitrary radii,
in $\reals^3$. The incidence graph $G(P,S)$ can be decomposed as
\begin{equation} \label{gps}
G(P,S) = G_0(P,S) \cup \bigcup_i (P_i\times S_i) ,
\end{equation}
such that, for each $i$, $P_i\subseteq P$ and $S_i\subseteq S$, and
\begin{align} \label{eq:main}
|G_0(P,S)| & = O\left( m^{2/3}n^{2/3} + m^{1/2}n^{7/8}\log^\beta(m^4/n) + m + n \right) , \\
\sum_i |P_i| & = O(m) , \quad\text{and}\quad \sum_i |S_i| = O(n)
\nonumber ,
\end{align}
where the constant exponent $\beta$ and the constants of
proportionality depend on the degree $D$ of $V$. Moreover, for each
$i$ there exists a circle $c_i\subset V$, such that $P_i = P\cap
c_i$ and $S_i$ is the set of spheres in $S$ that contain $c_i$.
\end{theorem}

\smallskip

\noindent{\bf Remark.} Apart from the improved bound on
$|G_0(P,S)|$, our bound is stronger, when compared to earlier works,
also in that it does not impose any restrictions on $G(P,S)$ (like
not containing a complete bipartite graph of some fixed size), and
gives a precise representation of graphs $G(P,S)$ that do have ``too
many'' incidences. (An earlier attempt at characterizing such large
graphs is given in Apfelbaum and Sharir~\cite{ApS} for the case of
planes (and hyperplanes in higher dimensions). Although it caters to
the general case (not requiring the points to lie on a surface), it
is much weaker than our representation.)

\paragraph{The cases where $V$ is (or contains) a plane or a sphere.}
We have explicitly ruled out these cases in our assumptions, because
the situation in these cases is different. These cases are treated
in Theorem~\ref{th:incgen} below. In these cases, each sphere
intersects $V$ in a circle, and the problem boils down to one
involving incidences between points and circles in the plane, or on
the sphere, except that the circles can have (potentially large)
multiplicities. Problems of this sort (without multiplicity of the
circles) have already been tackled in \cite{ANP+,ArS,MT}, and a
suitable extension of the analysis in these papers can also handle
multiplicities in a rather straightforward manner.

\paragraph{Applications.}
We begin by presenting two applications of Theorem~\ref{main}.
Actually, in the second application and in the first part of the
first one, we get better bounds, because the spheres that arise
there have a more constrained structure. First, we obtain the
following lower bounds on the number of distinct distances involving
points on a surface in $\reals^3$.
\begin{theorem} \label{dd3}
(a) Let $P$ be a set of $n$ points on an algebraic surface $V$ of
constant degree $D$ in $\reals^3$, which has no linear or spherical
components. Then the number of distinct distances determined by $P$
is $\Omega(n^{7/9}/\log^{\beta_1}n)$, where the constant exponent
$\beta_1$ and the constant of proportionality depend on
the degree $D$ of $V$. \\

\noindent (b) Let $P_1$ be a set of $m$ points on a surface $V$ as
in (a), and let $P_2$ be a set of $n$ arbitrary points in
$\reals^3$. Then the number of distinct distances determined by
pairs of points in $P_1\times P_2$ is
$$
\Omega \left( \min\left\{ m^{4/7}n^{1/7}/\log^{\beta_2}(m^4/n),\;
m,\; n \right\} \right) ,
$$
where the constant exponent $\beta_2$ and the constant of
proportionality depend on $D$.
\end{theorem}

While we believe that the bounds in the theorem are not tight, we
note that the bound in (a) is significantly larger than the
conjectured best-possible lower bound $\Omega(n^{2/3})$ for
arbitrary point sets in $\reals^3$, and so is the (somewhat weaker)
bound in (b), for a suitable range of values $m,n$ (including the
case $m=n$).

As another application, we bound the number of unit (or repeated)
distances involving points on a surface $V$, as above.
\begin{theorem} \label{und}
(a) Let $P$ be a set of $n$ points on some algebraic surface $V$ of
constant degree $D$ in $\reals^3$. Then $P$ determines $O(n^{4/3})$
unit distances,
where the constant of proportionality depends on the degree $D$ of $V$. \\

\noindent (b) Let $P_1$ be a set of $m$ points on a surface $V$ as
in (a), and let $P_2$ be a set of $n$ arbitrary points in
$\reals^3$. Then the number of unit distances determined by pairs of
points in $P_1\times P_2$ is
$$
O\left( m^{2/3}n^{2/3} + m^{6/11}n^{9/11} \log^{\beta_3} (m^3/n) + m
+ n \right) ,
$$
where the constant exponent $\beta_3$ and the constant of
proportionality depend on $D$.
\end{theorem}
The bound in (a) matches the best known upper bound for points in
the plane or on a sphere (see Brass, Moser, and Pach~\cite{BMP} for
a review of known results), and is in fact worst-case tight, since a
matching lower bound $\Omega(n^{4/3})$ is known for points on a
sphere with radius $1/\sqrt{2}$. The bound in (b) (say, for $m=n$)
is ``in between'' the best known general upper bound of $O(n^{3/2})$
for any set of $n$ points in $\reals^3$ (see~\cite{KMSS,Za}) and the
bound in (a).

Another interesting application of Theorem~\ref{main} is the
following general point-sphere incidence bound in three dimensions.
It improves the bound in Apfelbaum and Sharir~\cite{ApS2}, and is
more general, since it does not assume the spheres to be
non-degenerate, as is the case in \cite{ApS2}. This theorem is
analogous to the works of Brass and Knauer~\cite{BK} and Apfelbaum
and Sharir~\cite{ApS}, who studied incidences between points and planes 
(instead of spheres) in $\reals^3$ (and hyperplanes in higher dimensions).

\begin{theorem} \label{th:incgen}
Let $P$ be a set of $m$ points in $\reals^3$, and let $S$ be a set
of $n$ spheres, of arbitrary radii, in $\reals^3$. The incidence
graph $G(P,S)$ can be decomposed as
\begin{equation} \label{gpsc}
G(P,S) = G_0(P,S) \cup \bigcup_i (P_i\times S_i) ,
\end{equation}
such that, for each $i$, $P_i\subseteq P$ and $S_i\subseteq S$, and
\begin{align} \label{eq:mainc}
|G_0(P,S)| & = O\left( m^{8/11+\eps}n^{9/11} + m + n \right) , \\
\sum_i \big( |P_i| + |S_i| \big) & = O\left( m^{8/11+\eps}n^{9/11} + m + n \right) \nonumber ,
\end{align}
for any $\eps>0$, where the constant of proportionality depends on $\eps$.
Moreover, for each $i$ there exists a circle $c_i$, such that $P_i = P\cap c_i$ 
and $S_i$ is the set of spheres in $S$ that contain $c_i$.
\end{theorem}

\paragraph{The technique.}
Our approach continues the recent methodology of applying tools from
algebra and algebraic geometry to problems in combinatorial (and
computational) geometry, pioneered by Guth and Katz's
works~\cite{GK,GK2}. The main tool in this methodology is the
\emph{polynomial partitioning technique}, which yields a
divide-and-conquer mechanism via space decomposition, which in many
instances is a more effective tool than more traditional space
decomposition techniques (such as cuttings and simplicial
partitions; see, e.g., \cite{Chaz}). Interestingly though, while we
do use algebraic techniques, a major part of the analysis, involving
decomposition in a dual four-dimendsional space, goes back to the
source, and applies a standard cutting-based decomposition,
exploiting the fact that the objects that arise in this duality are
points and hyperplanes. This requires a somewhat more careful
analysis, but results in slightly improved bounds.

\smallskip

Theorem~\ref{main} will be proved in Section~\ref{sec:main},
Theorems~\ref{dd3} and~\ref{und} in Section~\ref{sec:dd3}, and
Theorem~\ref{th:incgen} in Section~\ref{sec:inc}.

\section{Proof of Theorem~\ref{main}} \label{sec:main}

Let $P$, $V$, $S$, $m$, and $n$ be as above. We first restrict the
analysis to the case where $V$ is irreducible. This can be done,
without loss of generality, by decomposing $V$ into its irreducible
components, assign each point of $P$ to each component that contains
it, and assign the spheres of $S$ to all the components. This
decomposes the problem into at most $D$ (as a matter of fact, at
most $D/2$) subproblems, each involving an irreducible surface, and
it thus follows that the original incidence count is at most $D/2$
times the bound for the irreducible case. In the remainder of this
section we thus assume that $V$ is irreducible.

To obtain the bound (\ref{eq:main}) on $I(P,S)$, we first derive a
weaker bound, and then improve it via a suitable decomposition of
dual space, similar to the way it has been done for circles in
\cite{ANP+,ArS}, and in more generality in \cite{SZ}, and also 
resembles the handling of the simpler case of points and planes in 
the companion paper \cite{SS16}.

\paragraph{A basic weak bound.}
By Sharir, Sheffer, and Zahl~\cite[Lemma 3.2]{SSZ}, except for at
most two \emph{``popular''} points, each point $p \in V$ is incident
to at most $44D^2=O(1)$ circles that are fully contained in $V$;
this follows since $V$ is neither a sphere nor a plane.

The number of incidences between the popular points, if any of them
is in $P$, and the spheres of $S$ is at most $2n$, so in what
follows we ignore these points and assume that $P$ does not contain
a popular point.

Let $C$ denote the set of circles that are fully contained in $V$,
contain at least one point of $P$, and are contained in at least one
sphere of $S$. For each circle $c\in C$ we form the bipartite
subgraph $P_c \times S_c$ of $G(P,S)$, where $P_c = P\cap c$ and
$S_c$ is the set of all the spheres of $S$ that contain $c$.

The preceding property therefore implies that $\sum_{{c}}
|P_{c}|=O(m)$. As $V$ is irreducible and non-spherical, it does not
contain any of the spheres in $S$. Thus, for each $s\in S$, the
intersection $s\cap V$ is an algebraic curve of degree at most $2D$
(as follows, e.g., from the generalized version of B\'ezout's
theorem~\cite{Fu84}), and can therefore contain at most $D=O(1)$
circles of $C$. This implies that $\sum_{{c}} |S_{c}|=O(n)$.

To recap, we have obtained a collection of complete bipartite graphs
$P_c\times S_c$, so that $\bigcup_{c} \left(P_{c}\times
S_{c}\right)$ is a portion of $G(P,S)$, $\sum_{c} |P_{c}| = O(m)$,
and $\sum_{c} |S_{c}| = O(n)$.

An interesting special case is when $V$ is \emph{ruled} by circles.
That is, each point $p\in V$ is incident to a circle that is fully
contained in $V$. Actually, as follows from a generalization of the
Cayley--Salmon theorem, established by Guth and Zahl~\cite{GZ}, an
irreducible surface of degree $D$ that is not ruled by circles can
fully contain at most $c D^2$ circles, for some absolute constant
$c$. That is, if $V$ is not ruled by circles, we also get a bound of
$O(D^2)=O(1)$ on the \emph{number} of complete bipartite graphs in
the decomposition~(\ref{gps}).

Surfaces ruled by circles contain infinitely many circles, but only
finitely many of them will yield nonempty bipartite graphs
$P_c\times S_c$. Informally, this means that we might get more
complete bipartite subgraphs in $G(P,S)$, but each with a smaller
number of edges; the linear bounds on the total size of their vertex
sets continue to hold.

For each $s\in S$, put $\gamma_s := (s\cap V) \setminus \bigcup C$.
As already observed, $V$ does not fully contain any sphere of $S$,
so each $\gamma_s$ is at most one-dimensional. By construction, it
does not contain any circle, and it might also be empty (for this or
for other reasons). Note that if $s\cap V$ does contain a circle
$c$, then ${c}\in C$, and the incidences between $s$ and the points
of $P$ on $c$ are all already recorded in $P_{c}\times S_{c}$
(assuming, of course, that $P_c = c\cap P\ne\emptyset$; otherwise,
removing $c$ incurs no loss of incidences). Finally, we ignore the
isolated points of $\gamma_s$. The number of
such points on a sphere $s$ is $O(1)$,\footnote{%
  The projection of $s\cap V$ onto a generic plane is a plane algebraic
  curve of degree at most $2D$~\cite{Harris}, and isolated points are
  projected onto isolated points. By Harnack's curve theorem~\cite{Har},
  the number of isolated points is thus $O(D^2) = O(1)$.}
so the number of incidences $(p,s)$, where $p$ is an isolated point
of $\gamma_s$, is at most $O(n)$.
We thus obtain the decomposition in (\ref{gps}), by letting
$G_0(P,S)$ denote the remaining portion of $G(P,S)$, after pruning
away the complete bipartite graphs $P_{c}\times S_{c}$. We further
remove from $G_0(P,S)$ all the $O(n)$ incidences involving isolated
points on their incident spheres, continue to denote the resulting
subgraph as $G_0(P,S)$, and put $I_0(P,S) = |G_0(P,S)|$.

Let $\Gamma$ denote the set of the $n$ curves $\gamma_s$, for $s\in
S$. The curves of $\Gamma$ are (spherical) algebraic curves of
degree at most $2D$ (e.g., see~\cite{Fu84}), and any pair of curves
$\gamma_s$, $\gamma_{s'}\in \Gamma$ intersect in at most $2D=O(1)$
points. Indeed, any of these points is an intersection point of $V$
with the circle $c=s\cap s'$; if $c$ is fully contained in $V$ then,
by construction, it has been removed from both curves, and if $c$ is
not contained in $V$, it can intersect it in at most $2D$ points.

Note that $I_0(P,S)$ is equal to the number of incidences
$I(P,\Gamma)$ between the points of $P$ and the curves of $\Gamma$.
To bound the latter quantity we proceed as follows.

We slightly tilt the coordinate frame to make it generic, and then
project the curves of $\Gamma$ onto the $xy$-plane. A suitable
choice of the tilting guarantees that (i) no pair of intersection
points or points of $P$ project to the same point, (ii) if $p$ is
not incident to $\gamma_s$ then the projections of $p$ and of
$\gamma_s$ remain non-incident, and (iii) no pair of curves in
$\Gamma$ have overlapping projections. In addition, by construction,
no curve of $\Gamma$ contains any (vertical) segment. Let $P^*$ and
$\Gamma^*$ denote, respectively, the set of projected points and the
set of projected curves; the latter is a set of $n$ plane algebraic
curves of constant maximum degree $2D$ (see, e.g.,
Harris~\cite{Harris} for the fact that projections do not increase
the degree). Moreover, $I(P,\Gamma)$ is equal to the number  $I(P^*,
\Gamma^*)$ of incidences between $P^*$ and $\Gamma^*$.

By the recent result of Sharir and Zahl~\cite{SZ} (see
Theorem~\ref{incPtCu} below), applied to $\Gamma^*$, the curves of
$\Gamma^*$ can be cut into $O(n^{3/2}\log^\kappa n)$ connected
Jordan subarcs, where the constant exponent $\kappa$ and the
constant of proportionality depend on $D$, so that each pair of
subarcs intersect at most once; the new arcs thus form a collection
of \emph{pseudo-segments}.

We can now apply the crossing-lemma technique of
Sz\'ekely~\cite{Sze}, exactly as was done in Sharir and
Zahl~\cite{SZ}. Since the resulting subarcs form a collection of
pseudo-segments, and the number of their intersections is $O(n^2)$,
Sz\'ekely's analysis yields the bound
$$
I_0(P,S) = I(P,\Gamma) = I(P^*,\Gamma^*) = O\left( m^{2/3}n^{2/3} +
m + n^{3/2}\log^\kappa n \right).
$$

Adding incidences recorded in the complete bipartite decomposition,
as constructed above, and the $O(n)$ incidences with isolated
points, we get our initial (weak) bound.
\begin{equation} \label{weakinc}
I(P,S) = O\left( m^{2/3}n^{2/3} + m + n^{3/2}\log^\kappa n +
\sum_{c} |P_{c}| \cdot |S_{c}| \right) ,
\end{equation}
where $\bigcup_{c} \left(P_{c}\times S_{c}\right)$ is contained in
the incidence graph $G(P,S)$, $\sum_{c} |P_{c}| = O(m)$, and
$\sum_{c} |S_{c}| = O(n)$.

\paragraph{The case $m=O(n^{1/4})$.}
Before proceeding to improve the bound in (\ref{weakinc}), we first
dispose of the case $m=O(n^{1/4})$. As above, we first remove all
the complete bipartite graphs $P_{c}\times S_{c}$, for ${c}\in C$,
from $G(P,S)$.
We then proceed to estimate $I_0(P,S)$ as follows. We call a
sphere $s\in S$ \emph{strongly degenerate} (or \emph{degenerate}\footnote{%
  This is much more restrictive than the notion of $\eta$-degeneracy mentioned earlier.}
for short) if all the points of $P\cap s$ are cocircular. We claim
that the number of incidences between the points of $P$ and the
non-degenerate spheres is $O(m^4+n)=O(n)$. Indeed, first discard the
spheres $s\in S$ containing at most three points of $P$, losing at
most $O(n)$ incidences. For an incidence between a point $p\in P$
and a surviving non-degenerate sphere $s\in S$, there exist (at
least) three distinct points $q,q',q'' \in (P\setminus \{p\})\cap s$
such that $p,q,q',q''$ do not all lie on a common circle; this
follows since $s$ is non-degenerate and $|s\cap P|\ge 4$, so $s$
contains at least three points that are not all cocircular with $p$.
The ordered quadruple $(p,q,q',q'')$ therefore \emph{uniquely}
accounts for the incidence between $p$ and $s$, and there are
$O(m^4)=O(n)$ such quadruples.

To bound the number of incidences between the points of $P$ and the
degenerate spheres of $S$, fix a degenerate sphere $s\in S$, and
assume that $m_s:=|s\cap P| \ge 2D+1$; the overall number of
incidences on the other spheres is at most $2Dn = O(n)$. By
assumption, all points of $s\cap P$ lie on a common circle $c$.
Since $c$ contains at least $2D+1$ points of $P$, it must be
contained in $V$, so the incidences involving $s$ are all recorded
in the complete bipartite graph $P_{c} \times S_{c}$. In other
words, we have shown, for $m = O(n^{1/4})$,
\begin{equation} \label{smallm}
I(P,S) = O\left( n + \sum_{c} |P_{c}| \cdot |S_{c}| \right) ,
\end{equation}
where, as above, $\bigcup_{c} \left(P_{c}\times S_{c}\right)$ is
contained in the incidence graph $G(P,S)$, and $\sum_{c} |P_{c}| =
O(m)$ and $\sum_{c} |S_{c}| = O(n)$.

\smallskip

\noindent{\bf Remark.} It seems likely that, with some care, the bound in (\ref{smallm})
could also be obtained by the technique of Fox et al.~\cite[Corollary 2.3]{FPSSZ} (see also \cite{SZ}).

\paragraph{Improving the bound.}
As in the analysis of incidences between points and circles (or
pseudo-circles) in \cite{ANP+,ArS}, and the more general analysis in
\cite{SZ}, the first two terms in (\ref{weakinc}) dominate when $m =
\Omega(n^{5/4}\log^{3\kappa/2} n)$. When $m$ is smaller, the third
term, which is independent of $m$, is the one that dominates, and we
then sharpen it as follows (a similar general approach is also used
in \cite{ANP+,ArS,SS16,SZ}).

We apply the standard lifting transform to 4-space, which maps each
point $(x,y,z)$ to the point $(x,y,z,x^2+y^2+z^2)$ on the paraboloid
$w=x^2+y^2+z^2$, and maps each sphere $(x-a)^2+(y-b)^2+(z-c)^2=r^2$
to the hyperplane $w=2ax+2by+2cz+(r^2-a^2-b^2-c^2)$. This lifting
preserves incidences: a point is incident to a sphere iff the lifted
point is incident to the lifted hyperplane. We next apply a standard
duality in 4-space that maps points to hyperplanes and vice versa,
and preserves point-hyperplane incidences. We denote by $a^*$ the
lifted and dualized image (hyperplane or point) of an object $a$
(point or sphere); to simplify the terminology, we call $a^*$ simply
the \emph{dual image} of $a$. 

We thus get a set $S^*$ of $n$ points,
and a set $P^*$ of $m$ hyperplanes in 4-space. We choose a parameter
$r$, to be fixed later, and construct a \emph{$(1/r)$-cutting} $\Xi$
for $P^*$ (see, e.g., \cite{Chaz}), which partitions $\reals^4$ into
$O(r^4)$ simplices, each crossed by at most $m/r$ dual hyperplanes.

\paragraph{Incidences with boundary dual points.}
Let us first handle dual points that lie on the boundaries of the
simplices of the $(1/r)$-cutting $\Xi$ and the dual hyperplanes.
\smallskip

\noindent (i) Each dual point that lies in the relative interior of
a 3-face $\varphi$ of some simplex of $\Xi$ has one incidence with
the dual hyperplane that contains $\varphi$ (if any), for a total of
$O(n)$ such incidences. If such a point $s^*$ is incident to another
dual hyperplane $p^*$, then $p^*$ must cross each simplex that has
$\varphi$ as a face (there are one or two such simplices). We then
assign $s^*$ to one of these simplices, and the relevant incidences
will then be counted in the subproblem associated with that simplex;
see below for details.

\smallskip

\noindent (ii) Consider next incidences involving dual points that
lie on some 2-face $f$ of some simplices of $\Xi$. In primal
4-space, the 2-flat $\pi$ spanned by $f$ is mapped to a line $\ell$,
such that any dual hyperplane $p^*$ that fully contains $f$ (that
is, $\pi$) is mapped back to a point in primal 4-space that lies in
$\ell$. Intersecting $\ell$ with the paraboloid $w=x^2+y^2+z^2$ and
projecting down to the original 3-space, we get at most two original
points $p$ whose dual images $p^*$ can fully contain $f$. Hence, the
number of incidences that fall into this special case are at most
$\sum_f 2|S^*\cap f|$, over all 2-faces $f$ as above.

\smallskip

\noindent (iii) For $f$, $\pi$, and $\ell$ as in (ii), the dual
hyperplanes $p^*$ that do not fully contain $f$ must cross every
simplex that has $f$ as a face. Indeed, let $\sigma$ be such a
simplex. Since $p^*$ meets $f$, it intersects the closure
$\bar{\sigma}$ of $\sigma$. Since $p^*$ does not cross $\sigma$, it
must be a supporting hyperplane to $\bar{\sigma}$. But such a
supporting hyperplane must meet $\bar{\sigma}$ in a full face of
some dimension. Hence, its intersection with $f$ must be at some
subface of $f$, contrary to assumption. Therefore, as in step (i),
we assign each dual point $s^* \in f$ to one of these adjacent
simplices, and the relevant incidences will then be counted in the
subproblem associated with that simplex.

Getting back to the bound $\sum_f 2|S^*\cap f|$ in (ii) for the
incidences with hyperplanes that fully contain $f$. Dual points
$s^*$ that lie on exactly one 2-face $f$ contribute a total of
$O(n)$ to the sum. Consider then a dual point $s^*$ that lies on two
(or more) 2-faces $f_1$, $f_2$. If $f_1$, $f_2$ are co-planar, we
can ignore one of them, because the two hyperplanes containing $f_1$
and the two containing $f_2$ are the same (see the argument in step
(ii) above, which only depends on the 2-flat supporting $f_1$ (which
is the same as the one supporting $f_2$), and not on $f_1$ (or
$f_2$) itself), so we still have at most two incidences involving
$s^*$. If $f_1$ and $f_2$ are not coplanar, then one of the
hyperplanes containing $f_2$ must cross $f_1$, so the two incidences
of $s^*$ within $f_2$ can be charged to this crossing incidence,
which is handled as above, by assigning $s^*$ into one of the
simplices bounded by $f_2$. It is easily checked that the same
argument applies when $s^*$ lies on more than two 2-faces: every
additional 2-face (which is not coplanar with $f_1$) will be
contained in a dual hyperplane that crosses $f_1$ (or else not
contribute any new incidence), so the corresponding incidences can
be charged to a suitable crossing incidence, as above. To summarize,
the sum $\sum_f 2|S^*\cap f|$ is $O(n)$ plus an excess that will be
handled within the simplices of $\Xi$.

\smallskip

\noindent (iv) Consider next incidences involving dual points on
some edge $e$ of some simplices of $\Xi$. In primal 4-space, the
line $\ell$ spanned by $e$ is mapped to a 2-flat $\pi$, such that
any dual hyperplane $p^*$ that fully contains $e$ (that is, $\ell$)
is mapped back to a point in primal 4-space that lies in $\pi$.
Intersecting $\pi$ with the paraboloid $w=x^2+y^2+z^2$ and
projecting down to the original 3-space, we conclude that the
original point $p$ that is mapped to $p^*$, for any $p^*$ as above,
lies in the intersection of $V$ with a circle $c$. Similarly, any
dual point $s^*$ that lies in $e$ (and thus in $\ell$) is mapped in
primal 4-space to a hyperplane that contains $\pi$. Intersecting
that hyperplane with the paraboloid, and projecting down to the
original 3-space, we obtain a sphere $s$ that fully contains $c$. We
may assume that $c$ is not fully contained in $V$, because
otherwise, incidences between points on $c$ and spheres that contain
$c$ are already recorded in the complete bipartite graph $P_c\times
S_c$ that we have removed from $G(P,S)$. But then $|P\cap c| \le 2D
= O(1)$. That is, at most $2D$ dual hyperplanes fully contain $e$,
yielding at most $2D|S^*\cap e|$ incidences, for a total, over all
edges $e$, of $\sum_e 2D|S^*\cap e|$ incidences.

\smallskip

\noindent (v) Arguing as in step (iii), dual hyperplanes $p^*$ that
cross $e$ but do not fully contain it must cross every simplex that
has $e$ as an edge. Hence, an incidence of such a dual hyperplane
$p^*$ with a dual point on $e$ can be charged to a crossing of some
adjacent simplex by $p^*$, and any such hyperplane-simplex crossing
can be charged only $O(1)$ times---at most once for each edge of the
simplex being crossed by $p^*$. In total, we get a total of
$O(r^4\cdot(m/r))=O(mr^3)$ such incidences.

Again, in the sum $\sum_e 2D|S^*\cap e|$, obtained in step (iv),
dual points $s^*$ that lie on just one simplex edge $e$ contribute
at most $2Dn=O(n)$ incidences. Consider then a point $s^*$ that lies
on more than one edge, say on edges $e_1$ and $e_2$. In the primal
space, they correspond to distinct circles $c_1$, $c_2$, both
contained in $s$, and any point $p$ on $c_2\setminus c_1$
corresponds to a dual hyperplane that (contains $e_2$ and) crosses
$e_1$, so the incidences with the dual hyperplanes that contain
$e_2$ can be charged to crossing incidences involving $e_1$, as
above. This also works for any number of edges containing $s^*$.
Hence, $\sum_e 2D|S^*\cap e|$ is $O(n)$ plus a term proportional to
the number of ``crossing incidences'', which, as has just been
argued, is $O(mr^3)$.

\smallskip

\noindent (vi) Consider finally incidences that involve dual points
that are vertices of some simplices of $\Xi$. If such a vertex $s^*$
does not lie in the relative interior of any higher-dimensional face
of any other simplex, that is, all the simplices adjacent to $s^*$
have $s^*$ as a vertex, then any incidence between a dual hyperplane
$p^*$ and $s^*$ can be charged to the crossing of $p^*$ with some
simplex $\sigma$ of $\Xi$ that is adjacent to $s^*$. It follows,
arguing as in (v), that the number of incidences of this kind is at
most $O(r^4)\cdot m/r = O(mr^3)$. On the other hand, if $s^*$ lies
in the relative interior of some higher-dimensional face $f$ of some
other simplex, we handle the incidence between $s^*$ and any
hyperplane $p^*$ as in steps (i)--(v) above.

To recap, ignoring incidences that have been assigned to the
subproblems within the simplices of $\Xi$, as well as incidences
that have been recorded in the complete bipartite graphs $P_c\times
S_c$, we have accumulated in this step only $O(n+mr^3)$ incidences.

\paragraph{Incidences within the simplices of $\Xi$.}
We now proceed to consider incidences within the simplices of $\Xi$.
For each simplex $\sigma$ of $\Xi$, let $n_\sigma$ denote the number
of points of $S^*$ in the interior of $\sigma$, including the points
that have been assigned to $\sigma$ from its boundary as above. We
bound, for each simplex $\sigma$ of $\Xi$, the number of incidences
between the $n_\sigma$ dual points in its interior and the at most
$m/r$ dual hyperplanes that cross $\sigma$. By duplicating simplices
$\sigma$ for which $n_\sigma > n/r^4$, so that in each copy we take
at most $n/r^4$ of these points (but retain all crossing
hyperplanes), we obtain a collection of $O(r^4)$ simplices, each of
which is crossed by at most $m/r$ dual hyperplanes and contains at
most $n/r^4$ dual points; we denote the actual number of these
hyperplanes and points as $m_\sigma$ and $n_\sigma$ (the latter
notation is slightly abused, as it now refers only to a single copy
(subproblem) of $\sigma$), respectively, for each simplex $\sigma$.

For each cell $\sigma$, we apply the bound (\ref{weakinc}) to the
subset $P^{(\sigma)}$ of the points of $P$ whose dual hyperplanes
cross $\sigma$, and to the subset $S^{(\sigma)}$ of the spheres
whose dual points lie in $\sigma$, and note that the case $m/r =
O((n/r^4)^{1/4})$ does not arise, because then we would also have
$m=O(n^{1/4})$, and then we would have used instead the bound
(\ref{smallm}), avoiding the partitioning altogether. That is, we
get, for each $\sigma$,
$$
I(P^{(\sigma)},S^{(\sigma)}) = O\left( m_\sigma^{2/3}n_\sigma^{2/3}
+ m_\sigma + n_\sigma^{3/2}\log^\kappa n_\sigma
 + \sum_c |P^{(\sigma)}_c| \cdot |S^{(\sigma)}_c| \right) ,
$$
for a suitable complete bipartite decomposition $\bigcup_c \left(
P^{(\sigma)}_c \times S^{(\sigma)}_c \right)$.

We sum these bounds, over the simplices $\sigma$ of $\Xi$. We note
that the same circle $c$ may arise in many complete bipartite graphs
$P^{(\sigma)}_c \times S^{(\sigma)}_c$, but (i) all these graphs are
contained in $P_c \times S_c$, and (ii) they are edge disjoint,
because each dual point $s^*$ lies in (or is a boundary point which
is assigned to) at most one simplex. This allows us to replace all
the partial subgraphs $P^{(\sigma)}_c \times S^{(\sigma)}_c$ by the
single graph $P_c \times S_c$, for each circle $c$ (contained in
$V$). We thus get

\begin{align} \label{incfinal}
I(P,S) & = O\left( \sum_\sigma \Big( m_\sigma^{2/3}n_\sigma^{2/3} +
 m_\sigma + n_\sigma^{3/2}\log^\kappa n_\sigma \Big)
 + mr^3 + n + \sum_c |P_c| \cdot |S_c| \right) \nonumber \\
& = O\left( r^4 \Big( (m/r)^{2/3}(n/r^4)^{2/3} +
(n/r^4)^{3/2}\log^\kappa (n/r^4) \Big)
 + mr^3 + n + \sum_c |P_c| \cdot |S_c| \right) \nonumber \\
& = O\left( m^{2/3}n^{2/3}r^{2/3} + \frac{n^{3/2}}{r^{2}}\log^\kappa
(n/r^4)
 + mr^3 + n + \sum_c |P_c| \cdot |S_c| \right) .
\end{align}

We now choose ${\displaystyle r = \frac{n^{5/16}\log^{3\kappa/8}
(m^4/n)}{m^{1/4}}}$, to equalize (asymptotically) the first two
terms in the bound (\ref{incfinal}), which then become
$O(m^{1/2}n^{7/8}\log^{\kappa/4}(m^4/n))$. The third term becomes
$mr^3 = m^{1/4}n^{15/16}\log^{9\kappa/8} (m^4/n)$, which is
dominated by the preceding bound for $m=\Omega(n^{1/4})$, as is
easily checked. As already noted, the complementary case
$m=O(n^{1/4})$ has been handled by (\ref{smallm}), and the case $m =
\Omega(n^{5/4}\log^{3\kappa/2}n)$ is handled simply by
(\ref{weakinc}) (now without the term $n^{3/2}\log^\kappa n$ as it
is subsumed by the other term).

This completes the proof of Theorem~\ref{main}. $\Box$

\smallskip

\noindent{\bf Remark.} In retrospect, once we have reduced the
problem to that of bounding the number $I(P^*,\Gamma^*$) on
incidences between the projected points and curves on the
$xy$-plane, we could have applied, as a black-box, the analysis of
Sharir and Zahl~\cite{SZ}, and get a slightly weaker bound with an
additional (arbitrarily small) $\eps$ in the exponents.

As this remark will be significant in the proofs of some of the
forthcoming applications, we rephrase here the result of~\cite{SZ},
with the notation used in the proof of Theorem~\ref{main}, for the
convenience of the reader.

\begin{theorem}[Sharir and Zahl~\cite{SZ}]\label{incPtCu}
Let $\Gamma^*$ be a set of $n$ algebraic plane curves that belong to
an $s$-dimensional family $F$ of curves of maximum constant degree
$D$, no two of which share a common irreducible component, and let
$P^*$ be a set of $m$ points in the plane. Then, for any $\eps>0$,
the number $I(P^*,\Gamma^*)$ of incidences between the points of
$P^*$ and the curves of $\Gamma^*$ satisfies
\begin{equation*}
 I(P^*,\Gamma^*) = O\Big(m^{\frac{2s}{5s-4}} n^{\frac{5s-6}{5s-4}+\eps} + m^{2/3}n^{2/3} + m + n\Big) ,
\end{equation*}
where the constant of proportionality depends on $\eps$, $s$, $D$,
and the complexity of the family $F$.
\end{theorem}
In this result, an \emph{$s$-dimensional family of curves} is a family $\C$
of algebraic curves (of constant maximum degree), so that each curve
in $\C$ can be represented as a point in some finite-dimensional
parametric space, and the set of these ``dual'' points is an
$s$-dimensional algebraic variety of constant degree (which is
referred to as the ``complexity'' of $F$).

In our case, $s=4$, since each curve of $\Gamma^*$ can be
represented by the four parameters that define the corresponding
sphere, and, using the fact that $V$ is of constant degree, it is
easy to verify that the assumptions of Theorem~\ref{incPtCu} hold in
this case. Substituting $s=4$ gives us our bound, except that the
polylogarithmic factor is replaced by the factor $n^\eps$. That is,
exploiting the fact that we are dealing here with spheres, so that
the dual representation involves points and \emph{hyperplanes},
allows us to obtain the finer bound in (\ref{eq:main}) (concretely,
using cuttings instead of the rather involved partitioning scheme of
\cite{MP}).

\section{Distinct and repeated distances in three dimensions}
\label{sec:dd3}

In this section we prove Theorems~\ref{dd3} and \ref{und}, the
applications of our main result to distinct and repeated distances 
in three dimensions. The theorems present four results, in each of which the 
problem is reduced to one involving incidences between spheres and points on a
surface. However, except for Theorem~\ref{dd3}(b), the spheres
that arise in the other cases are restricted, by requiring their
centers to lie on $V$ and / or to have a fixed radius. This makes
the spheres have only three or two degrees of freedom. The case of
two degrees of freedom (in Theorem~\ref{und}(a)) is the simplest,
and requires very little of the machinery developed here (see
below). The cases of three degrees of freedom (in
Theorem~\ref{dd3}(a) and Theorem~\ref{und}(b)) call for a dual
representation of the setup in three dimensions.

A rigorous analysis along this line is doable, and we will comment
on it later, but there are several technical issues that arise, and
a careful treatment of them will make the proofs longer and somewhat
more involved. As a compromise, we state the sharp bounds that would
result from the full analysis, but present simpler ``black-box''
proofs that exploit the machinery in \cite{SZ} and yield slightly
inferior bounds.

\noindent{\bf Proof of Theorem~\ref{dd3}.} We will first establish
the more general bound in (b); handling (a) requires a somewhat
different analysis.

\smallskip

\noindent{\bf (b)} Let $t$ denote the number of distinct distances
in $P_1\times P_2$. For each $q\in P_2$, draw $t$ spheres centered
at $q$ and having as radii the $t$ distinct distances. We get a
collection $S$ of $nt$ spheres, a set $P_1$ of $m$ points on $V$,
which we relabel as $P$, to simplify the notation, and exactly $mn$
incidences between the points of $P$ and the spheres of $S$.

In order to effectively apply the bound in Theorem~\ref{main}, we
first have to control the term $\sum_i |P_{i}|\cdot |S_{i}|$; that
is, we argue that most of the $mn$ incidences do not come from this
bound, unless $t=\Omega(n)$. Indeed, for each $i$, we have $|S_i|
\le 2t$; this is because all the spheres in $S_i$ pass through a
fixed circle $c$, so, up to multiplicity $2$, their radii are all
distinct. This implies that
$$
\sum_i |P_{i}|\cdot |S_{i}| \le 2t\sum_i |P_i| = O(mt) .
$$
If this would have accounted for more than, say, half the
incidences, we would get $t =\Omega(n)$, as claimed, and then the
bound in the theorem would follow. We may thus ignore this term, and
write
$$
mn = O\left( m^{2/3}(nt)^{2/3} + m^{1/2}(nt)^{7/8} \log^\beta(m^4/n)
+ m + nt \right) ,
$$
or
$$
t = \Omega \left( \min\left\{ m^{1/2}n^{1/2},\;
m^{4/7}n^{1/7}/\log^{8\beta/7}(m^4/n),\; m \right\} \right) ,
$$
as claimed.

\smallskip

\noindent{\bf (a)} Here we are in a more favorable situation,
because the spheres in $S$ have only three degrees of freedom, in
the sense that their centers lie on the surface $V$, so that, in
principle, we need only two parameters to specify the center and one
for the radius.

One possibility would be to adapt the analysis in the proof of
Theorem~\ref{main}, with the difference that the spheres are now
dualized to points in three dimensions, rather than four. As already
noted, this would raise several technical issues, which, albeit
minor, require careful analysis that would be too space-consuming. A
discussion of the issues that arise and the way to handle them to
get the sharper bound is given below.

Instead, we ``shortcut'' the analysis, and apply the improved
incidence bound of Sharir and Zahl~\cite{SZ}, stated in
Theorem~\ref{incPtCu}, with $s=3$. That is, we still represent each
curve $\gamma_s^*$ of $\Gamma^*$ by the parameters $(x,y,z,r)\in
\reals^4$ that define the corresponding sphere $s$ (where $(x,y,z)$
is its center and $r$ its radius), but now $(x,y,z)$ is constrained
to lie on $V$. It then easily follows that $\Gamma^*$ is a
three-dimensional family of curves (in the notation of
Theorem~\ref{incPtCu}).

We thus get the bound
$$
I(P^*,\Gamma^*) = O\left( |P^*|^{2/3}|\Gamma^*|^{2/3} +
|P^*|^{6/11}|\Gamma^*|^{9/11+\eps} +|P^*| + |\Gamma^*| \right) ,
$$
for any $\eps>0$. Arguing as in the proof of (b), we may ignore the
term $\sum_i |P_i|\cdot |S_i|$ in the bound on $I(P,S)$, which is
negligible unless $t=\Omega(n)$, and thus get the inequality
$$
n^2 = O\left( n^{2/3}(nt)^{2/3} + n^{6/11}(nt)^{9/11+\eps} + nt
\right) ,
$$
which yields $t=\Omega(n^{7/9-\eps})$, for any $\eps>0$, thereby
completing the proof of (the coarser version of) (a). $\Box$

\noindent{\bf Proof of Theorem~\ref{und}.}

Consider (a) first. Following the standard approach to problems
involving repeated distances, we draw a unit sphere $s_p$ around
each point $p\in P$, and seek an upper bound on the number of
incidences between these spheres and the points of $P$; this latter
number is exactly twice the number of unit distances determined by
$P$.

This instance of the problem has several major advantages over the
general analysis in Theorem~\ref{main}. First, in this case the
incidence graph $G(P,S)$ cannot contain $K_{3,3}$ as a subgraph,
eliminating altogether the complete bipartite graph decomposition in
(\ref{gps}) (or, rather, bounding the overall number of edges in
these subgraphs by $O(n)$).

More importantly, the family $\Gamma^*$ of curves ``almost'' has
only two degrees of freedom. To have two degrees of freedom, in the
sense of Pach and Sharir~\cite{PS}, it is required that, for any
pair of points $p^*,q^*\in P^*$, there are at most $O(1)$ curves of
$\Gamma^*$ passing through $p^*$ and $q^*$ (and that any pair of
curves of $\Gamma^*$ intersect in at most $O(1)$ points, a property
that we have already established).

To test for this property, fix a pair $p^*,q^*\in P^*$. By our
assumption that the coordinate frame is generic, there is a unique
pair $p,q\in P$ that project, respectively, to $p^*$ and $q^*$, and
any curve $\gamma_s^*\in\Gamma^*$ that passes through $p^*$ and
$q^*$ is the projection of a unique curve $\gamma_s\in\Gamma$ that
passes through $p$ and $q$. The corresponding sphere $s\in S$ is
then a unit sphere that passes through $p$ and $q$, so its center
must lie on a suitable circle $c_{pq}$ that is centered at
$\frac12(p+q)$ and is orthogonal to $\vec{pq}$. As is easily
checked, the circles $c_{pq}$ are all distinct.

If $c_{pq}$ is not fully contained in $V$, it meets it in at most
$2D$ points, implying that there are at most $2D=O(1)$ curves of
$\Gamma^*$ that pass through $p^*$ and $q^*$, as desired.

It remains to study pairs $p^*,q^*$ for which $c_{pq}$ is fully
contained in $V$. The number of curves that pass through $p^*$ and
$q^*$ is $|P\cap c_{pq}|$. By the result of Sharir, Sheffer, and
Zahl~\cite{SSZ}, already mentioned in the proof of
Theorem~\ref{main}, except for two popular points (which we may
assume, as above, not to belong to $P$), every point $p\in P$ is
incident to at most $44D^2=O(1)$ circles that are fully contained in
$V$. It follows that
$$
\sum_{p,q\in P} |P\cap c_{pq}| = O(n) .
$$
We can now apply Sz\'ekely's crossing lemma argument~\cite{Sze} to
$P^*$ and $\Gamma^*$. The edges in Sz\'ekely's graph have constant
multiplicity, except for those that connect pairs $p^*, q^*$ for
which $c_{pq}\subset V$. As just argued, the overall number of edges
of the latter kind is $O(n)$. Omitting these edges from the graph,
Sz\'ekely's argument applies to the remainder, and yields the bound
$O(n^{4/3})$ for the number of edges. Combining this bound with the
linear bound on the number of high-multiplicity edges, and the
additional linear bound on the size of the complete bipartite graphs
$P_i\times S_i$, as noted above, we get a total of $O(n^{4/3})$
incidences, and thus $O(n^{4/3})$ unit distances.

\medskip

\noindent We now consider (b). Again, we reduce the problem to that
of bounding the number of incidences between the $m$ points of
$P_1$, which lie on $V$, and the $n$ unit spheres centered at the
points of $P_2$. Here too the overall number of edges in the
complete bipartite graph decomposition is $O(m+n)$, so we can ignore
this part of the bound.

In this case, the curves of $\Gamma^*$ have three degrees of
freedom, or, in the terminology of Sharir and Zahl~\cite{SZ}, as
reviewed in Theorem~\ref{incPtCu}, $\Gamma$ is a three-dimensional
family of curves. Applying the same reasoning as in this preceding
proof, we conclude that the number of unit distances in this case is
$$
O\left( m^{2/3}n^{2/3} + m^{6/11}n^{9/11+\eps} + m + n \right) ,
$$
for any $\eps>0$. 
$\Box$

\paragraph{Improving the bounds.}
In the proofs of Theorem~\ref{dd3}(a) and Theorem~\ref{und}(b), we
want to dualize the problem in a way that exploits the fact that the
spheres of $S$ have only three degrees of freedom. We still map the
spheres to points in $\reals^4$ and the points to hyperplanes in
$\reals^4$, as above, but now the dual points $s^*$ all lie on a
three-dimensional algebraic variety $V^*\subset\reals^4$ of constant
degree; in Theorem~\ref{dd3}(a), $V^* = V\times\reals$, and in
Theorem~\ref{und}(b), $V^*$ is the paraboloid $x_4 =
x_1^2+x_2^2+x_3^2+1$. We construct a $(1/r)$-cutting for the
collection of the dual hyperplanes, but use a coarser (and simpler)
technique of drawing a random sample of $O(r\log r)$ hyperplanes,
construct their arrangement, and triangulate each cell into
simplices. We now use the generalized zone theorem of Aronov et
al.~\cite{APS}, to conclude that $V^*$ crosses only $O(r^3\log^4r)$
simplices, and we apply the weak bound only in these simplices.
There are various additional technical issues that have to be
handled, but, as already explained above, we skip over them, in the
interest of keeping the paper short. Working out all the details, we
get the slightly improved bounds, as asserted in the theorems.

\section{Proof of Theorem~\protect{\ref{th:incgen}}} \label{sec:inc}

The proof establishes the bound in (\ref{eq:mainc}), via induction,
adding a prespecified approximation parameter $\eps>0$ to the bound.
Concretely, we claim that, for any prespecified $\eps>0$ we can
write
$$
G(P,S) = G_0(P,S) \cup \bigcup_c (P_c\times S_c) ,
$$
so that
\begin{equation} \label{eq:8911}
|G_0(P,S)| \le A\left( m^{8/11+\eps}n^{9/11} + m + n \right) ,
\end{equation}
and
$$
\sum_c \big( |P_c| + |S_c| \big) = O\left( m^{8/11+\eps}n^{9/11} + m
+ n \right) ,
$$
where $A$ and the other constants of proportionality depend on $\eps$. 

The base cases are when $m\le n^{1/4}$ or when $m\le m_0$, for some
sufficiently large constant $m_0$ that will be set later.

Consider first the case $m\le n^{1/4}$. We adapt the argument for
this case given in the proof of Theorem~\ref{main}. It yields the
bound $I(P,S_1)=O(n)$ for the set $S_1$ of spheres of $S$ that are
not strongly degenerate. Each strongly degenerate sphere $s$ can be
replaced by the unique circle $c_s$ that contains all its incident
points. We thus get an incidence problem involving a set $P$ of $m$
points and a multiset $\C$ of $n$ circles in $\reals^3$.

Fix the threshold multiplicity $\mu_0=n^{1/4}$, and consider the set
$\C^-$ of all circles in $\C$ with multiplicity at most $\mu_0$.
Each circle $c\in\C^-$ with at most two points of $P$ on it
contributes at most $2\mu(c)$ incidences, where $\mu(c)$ denotes the
multiplicity of $c$. Summing these bounds over all such circles, we
get at most $2|S|=2n$ incidences. The number of incidences involving
circles containing at least three points of $P$ is $O(m^3\cdot\mu_0)
= O(n)$.

This leaves us with circles of multiplicity larger than $n^{1/4}$.
We represent the corresponding incident pairs as a union of complete
bipartite graphs $P_c\times S_c$, over all circles in $\C$ with
multiplicity larger than $n^{1/4}$. We clearly have $\sum_c |S_c|
\le n$, and $\sum_c |P_c|$ is simply the number of incidences
between the points of $P$ and the at most $n^{3/4}$ ``heavy''
circles, counted \emph{without} multiplicity. The same argument used
above gives the bound $O(m^3+n^{3/4}) = O(n^{3/4})$.

In summary, we have for this case
\begin{equation} \label{smallmc}
I(P,S) = O\left( n + \sum_{c} |P_{c}| \cdot |S_{c}| \right) ,
\end{equation}
where $\bigcup_{c} \left(P_{c}\times S_{c}\right)$ is contained in
the incidence graph $G(P,S)$, and $\sum_{c} |P_{c}| = O(n^{3/4})$
and $\sum_{c} |S_{c}| = O(n)$. That is,~(\ref{eq:8911}) holds in
this case.

The case $m\le m_0$ follows easily if we choose $A$ sufficiently
large. This holds for any choice of $m_0$; the value that we choose
is specified later.

Suppose then that (\ref{eq:8911}) holds for all sets $P'$, $S'$,
with $|P'|<m$, $|S'|<n$, and consider the case where the sets $P,S$
are of respective sizes $m,n$, and $m>n^{1/4}$ and $m>m_0$.

Before continuing, we also dispose of the case $m\ge n^3$. In this case
we consider the arrangement $\A(S)$ of the spheres in $S$. The complexity
of $\A(S)$ is $O(n^3)$. More precisely, this bound holds, and is 
asymptotically tight, for spheres in general position. In our case, 
$S$ is likely not to be in general position, and then the complexity of
$\A(S)$ might be smaller, because vertices and edges might be incident to
many spheres. Nevertheless, if we count each vertex and edge of $\A(S)$
with its multiplicity, we still get the upper complexity bound $O(n^3)$.

This means that the number of incidences with points that are either
vertices or lie on the (relatively open) faces of $\A(S)$ is $O(n^3)=O(m)$.
Incidences with points that lie on the (relatively open) edges of $\A(S)$ 
(note that each such edge is a portion of some circle) are recorded, as usual, 
by a complete bipartite graph decomposition $\bigcup_c (P_c\times S_c)$,
where, as just argued, we have $\sum_c |P_c| \le m$ and
$\sum_c |S_c| = O(n^3) = O(m)$. This implies that (\ref{eq:8911}) holds 
in this case, so, in what follows, we assume that $m\le n^3$.

\paragraph{Applying the polynomial partitioning technique.}
We fix a sufficiently large constant parameter $D \ll m^{1/3}$,
whose choice will be specified later, and apply the polynomial
partitioning technique of Guth and Katz~\cite{GK2}. We obtain a
polynomial $f\in\reals[x,y,z]$ of degree at most $D$, whose zero set
$Z(f)$ partitions 3-space into $O(D^3)$ (open) connected components
(cells), and each cell contains at most $O(m/D^3)$ points. By
duplicating cells if necessary, we may also assume that each cell is
crossed by at most $O(n/D)$ spheres of $S$; this duplication keeps
the number of cells $O(D^3)$ (because each sphere crosses only $O(D^2)$ cells). 
That is, we obtain at most $aD^3$
subproblems, for some absolute constant $a$, each associated with
some cell of the partition, so that, for each $i\le aD^3$, the
$i$-th subproblem involves a subset $P_i\subset P$ and a subset
$S_i\subset S$, such that $m_i:=|P_i|\le bm/D^3$ and $n_i:=|S_i|\le
bn/D$, for another absolute constant $b$.

Set $P_0:= P\cap Z(f)$ and $P'=P\setminus P_0$. We have
\begin{equation} \label{eqsplit}
I(P,S) = I(P_0,S) + I(P',S) .
\end{equation}
We first bound $I(P_0,S)$. Decompose $Z(f)$ into its $O(D)$
irreducible components, assign each point of $P_0$ to every
component that contains it, and assign the spheres of $S$ to all
components. We now fix a component, and bound the number of
incidences between the points and spheres assigned to that
component; $I(P_0,S)$ is at most $D$ times the bound that we get.

We may therefore assume that $Z(f)$ is irreducible. If $Z(f)$ is a
plane or a sphere, then for any sphere $s \in S$, the curve $s\cap
Z(f)$ is a circle; let $\C$ denote the multiset of these circles,
where each circle has multiplicity equal to the number of spheres
that contain it. Then $I(P_0,S)$ is the number of incidences between
the points of $P_0$ and the circles of $\C$, counted with
multiplicity.

We bound the number of incidences of this latter kind using the
incidence bound of Aronov et al.~\cite{AKS} for points and circles
in $\reals^3$. Fixing a threshold $\mu$, the number of incidences
involving circles with multiplicity at most $\mu$ (and counted with
their multiplicity) is easily seen to be
\begin{equation} \label{pcmult}
O\left( m^{2/3}n^{2/3}\mu^{1/3} +
m^{6/11}n^{9/11}\mu^{2/11}\log^{2/11} (m^3\mu/n) + m\mu + n \right)
.
\end{equation}
We now choose
$$
\mu = \min \Big\{ m^{2/11}n^{5/11},\; m,\; n^{9/11}/m^{3/11} \Big\}
.
$$
An easy, albeit a bit tedious, calculation shows that the bound in
(\ref{pcmult}) becomes $O(m^{8/11}n^{9/11} + n)$.

For circles $c$ with multiplicity larger than $\mu$, we record the
corresponding point-sphere incident pairs by a complete bipartite
graph decomposition $\bigcup_c (P_c\times S_c)$, where $c$ ranges over all
such ``heavy'' circles, and where $P_c = P_0\cap c$ and $S_c$ is the
set of all spheres that contain $c$. We clearly have $\sum_c |S_c| =
O(n)$ (each sphere can intersect $Z(f)$ in only one circle, except for
the unique sphere, if any, that coincides with $Z(f)$, which we may ignore),
and $\sum_c |P_c|$ is the number of incidences between the
points of $P_0$ and the heavy circles, counted \emph{without}
multiplicity. The number of these circles is at most $O(n/\mu)$.
Using the bound in \cite{AKS}, we get, as above,
$$
\sum_c |P_c| = O\left( \frac{m^{2/3}n^{2/3}}{\mu^{2/3}} +
\frac{m^{6/11}n^{9/11}}{\mu^{9/11}}\log^{2/11} (m^3\mu/n) + m +
\frac{n}{\mu} \right) .
$$
Since this is asymptotically the same as the bound in (\ref{pcmult})
divided by $\mu$, we simply (and pessimistically) upper bound this
by $O(m^{8/11}n^{9/11} + n)$.

Assume then that $Z(f)$ is neither a plane nor a sphere. Since
$\deg(Z(f))\le D$ is a constant, our main Theorem~\ref{main} implies
that
$$
I(P_0,S) = O\left( m^{2/3}n^{2/3} + m^{1/2}n^{7/8}\log^\beta(m^4/n)
+ m + n + \sum_c |P_c|\cdot |S_c| \right) ,
$$
where $\bigcup_c (P_c\times S_c) \subseteq G(P_0,S)$, and $\sum_c
|P_c| = O(m)$ and $\sum_c |S_c| = O(n)$. As is easily checked, the
first four terms are subsumed in~(\ref{eq:8911}), if we choose $A$
sufficiently large, and the term $\sum_c |P_c|\cdot |S_c|$ is added
to the complete bipartite graph decomposition that we collect.

Finally, we estimate
$$
I(P',S) = \sum_{i=1}^{aD^3} I(P_i,S_i) .
$$
By the induction hypothesis, we get
$$
I(P_i,S_i) \le A \left(m_i^{8/11+\eps}n_i^{9/11} + m_i + n_i +
\sum_c |P_{i,c}| \cdot |S_{i,c}| \right) ,
$$
for a suitable complete bipartite decomposition $\bigcup_c \left(
P_{i,c} \times S_{i,c} \right) \subseteq G(P_i,S_i)$.

When summing these bounds, we note that the same circle $c$ may
arise in many complete bipartite graphs, but, as already noted
earlier, (i) all these graphs are contained in $P_c \times S_c$, and
(ii) they are edge disjoint, because each point $p \in P'$ lies in
at most one cell, and even if this cell gets duplicated, the
relavant spheres are all distinct. This allows us to replace all the
partial subgraphs $P_{i,c} \times S_{i,c}$ by the single graph $P_c
\times S_c$, for each circle $c$ in the decomposition.

The sum of the other terms is
\begin{align*}
I_0(P',S) & \le A\cdot aD^3 \left( (bm/D^3)^{8/11+\eps} (bn/D)^{9/11} + (bm/D^3) + (bn/D) \right) \\
& = \frac{Aab^{17/11+\eps}}{D^{3\eps}} m^{8/11+\eps}n^{9/11} + Aabm
+ AabD^2n .
\end{align*}
We note that $m^{8/11+\eps}n^{9/11} \ge m^\eps\cdot m$ and
$m^{8/11+\eps}n^{9/11} \ge m^\eps\cdot n$ for for $n^{1/4}\le m\le
n^3$. We choose $D$ sufficiently large so that $D^{3\eps} \ge
4ab^{17/11+\eps}$, and then the bound is at most
$$
\left( \frac{A}{4} + \frac{Aab}{m^\eps} + \frac{AabD^2}{m^\eps}
\right) m^{8/11+\eps}n^{9/11} .
$$
Choosing $m_0$ sufficiently large, so that $m_0^\eps \ge 4abD^2$, we
ensure that, for $m\ge m_0$,
$$
\frac{A}{4} + \frac{Aab}{m^\eps} + \frac{AabD^2}{m^\eps} \le
\frac{A}{4} + \frac{A}{4} + \frac{A}{4} = \frac{3A}{4} .
$$
Adding the bounds for $I(P_0,S)$, and choosing $A$ sufficiently
large, we get that (\ref{eq:8911}) hold for $P$ and $S$. The
coorresponding bounds on $\sum_c |P_c|$ and $\sum_c |S_c|$ are
established by the same inductive analysis, and we omit the
straightforward details. This establishes the induction step and
thereby completes the proof. $\Box$

\section{Discussion}

In this paper we have made significant progress on a major incidence
problem involving points and spheres in three dimensions, for the
special case where the points lie on a constant-degree algebraic
surface. We have also obtained several applications of this result
to problems involving repeated and distinct distances in three
dimensions, and have extended the analysis to the case where the
points are arbitrary and are not required to lie on a
constant-degree surface; this latter extension improves the
bound derived in Apfelbaum and Sharir~\cite{ApS2}, and it is also
significantly more general, as it does not require the spheres to be
non-degenerate, as in \cite{ApS2}.

The study in this paper raises several interesting open problems.

%

\smallskip

\noindent{(i)} Our analysis suggests that if, instead of the family
of spheres, we take $S$ to be a \emph{$k$-dimensional family} of
constant-degree algebraic surfaces (in the terminology of \cite{SZ},
already mentioned above), and still assume the points to lie on a
constant-degree surface, we can then extend the analysis in
Theorem~\ref{main} to get an analogous bound for point-surface
incidences, depending on $k$, and resembling the one obtained in
\cite{SZ} for point-curve incidences in the plane.

It also seems likely that, as in Theorem~\ref{th:incgen}, the analysis
can be further extended to the case where the points do not have to
lie on a surface, and that the bound that it yields is
$$
O\left( m^{\frac{2s}{3s-1}+\eps}n^{\frac{3s-3}{3s-1}} + m + n +
\sum_\gamma |P_\gamma|\cdot |S_\gamma| \right) ,
$$
for any $\eps>0$, where $\bigcup_\gamma (P_\gamma\times S_\gamma)
\subseteq G(P,S)$, and where $\sum_\gamma |P_\gamma|$, $\sum_\gamma
|S_\gamma|$ are suitably bounded. Assuming that this extension can
be made rigorous, it would yield a significant generalization of
Zahl's result~\cite{Za}, where the leading term almost matches his
bound, but there are no restrictions that the incidence graph does
not contain a fixed-size complete bipartite subgraph, as assumed in
\cite{Za}.

\smallskip

\noindent{(ii)} A long-standing open problem is that of establishing
the lower bound of $\Omega(n^{2/3})$ for the number of distinct
distances determined by a set of $n$ points in $\reals^3$, without
assuming them to lie on a surface. The best known lower bound is due
to Solymosi and Vu~\cite{SoVu}. In the present study we have
obtained some partial results (with better lower bounds) for cases
where the points do lie on a surface. We hope that some of the ideas
used in this work could be applied for the general problem.

\smallskip

\noindent{(iii)} Another major long-standing open problem is that of
improving the upper bound $O(n^{3/2})$, established in
\cite{KMSS,Za}, on the number of unit distances determined by a set
of $n$ points in $\reals^3$, again without assuming them to lie on a
surface. It would be interesting to make progress on this problem.

\smallskip

\noindent{(iv)} Finally, it would be interesting to find additional
applications of the results of this paper. One direction to look at
is the analysis of repeated patterns in a point set, such as
congruent or similar simplices, which can sometimes be reduced to
point-sphere incidence problems; see \cite{AAPS,AgS}.

\end {document}